\documentclass[amstex,12pt, amssymb]{article}

\usepackage{mathtext}
\usepackage[cp1251]{inputenc}
\usepackage[T2A]{fontenc}
\usepackage[dvips]{graphicx}
\usepackage{amsmath}
\usepackage{amssymb}
\usepackage{amsxtra}
\usepackage{latexsym}
\usepackage{ifthen}

\newcounter{lemma}[section]

\newcounter{corollary}[section]

\newcounter{remark}[section]

\newcounter{theorem}[section]

\newcounter{proposition}[section]

\newcounter{example}

\numberwithin{equation}{section}

\begin{document}

\markboth{\centerline{O.~DOVHOPIATYI, E.~SEVOST'YANOV}}
{\centerline{ON COMPACT CLASSES...}}

\def\cc{\setcounter{equation}{0}
\setcounter{figure}{0}\setcounter{table}{0}}

\overfullrule=0pt

%\normalsize\large

\author{{OLEKSANDR DOVHOPIATYI, EVGENY SEVOST'YANOV}\\}

\title{
{\bf ON COMPACTNESS OF ONE CLASS OF SOLUTIONS OF THE DIRICHLET
PROBLEM}}

\date{\today}
\maketitle

%\large
\begin{abstract}
We consider the Dirichlet problem for the Beltrami equation in some
simply connected domain. We consider the class of all homeomorphic
solutions of such a problem with a normalization condition and
set-theoretic constraints on their complex characteristics. We have
proved the compactness of this class in terms of prime ends for an
arbitrary continuous function in the Dirichlet condition.
\end{abstract}


\section{Introduction}

Relatively recently, different authors have obtained theorems on the
compactness of the classes of solutions of the Beltrami equations
with different types of restrictions on their complex
characteristics, see, for example, \cite{Dyb$_1$}, \cite{LGR} and
\cite{L$_1$}--\cite{L$_2$}. In particular, the results of
\cite{L$_1$}--\cite{L$_2$} were used in the variational method,
which is reflected in full volume in the article \cite{LGR}. The
problem of the compactness of the classes of solutions of the
Dirichlet problem for this equation is quite close in relation to
the question of the compactness of solutions to the Beltrami
equation. The authors of this manuscript obtained some results on
this topic for the case when the area in which the Dirichlet problem
is studied is Jordanian, see, for example, \cite{DS} and \cite{SD}.
This article contains results on the compactness of the classes of
solutions to the Dirichlet problem for the case when the given
domain is simply connected (not necessarily Jordan), and the complex
characteristics of the solutions satisfy set-theoretic constraints.

\medskip
In what follows, a mapping $f:D\rightarrow{\Bbb C}$ is assumed to be
{\it sense-preserving,} moreover, we assume that $f$ has partial
derivatives almost everywhere. Put $f_{\overline{z}} = \left(f_x +
if_y\right)/2$ and $f_z = \left(f_x - if_y\right)/2.$ The {\it
complex dilatation} of $f$ at $z\in D$ is defined as follows:
$\mu(z)=\mu_f(z)=f_{\overline{z}}/f_z$ for $f_z\ne 0$ and $\mu(z)=0$
otherwise. The {\it maximal dilatation} of $f$ at $z$ is the
following function:
\begin{equation}\label{eq1}
K_{\mu}(z)=K_{\mu_f}(z)=\quad\frac{1+|\mu (z)|}{1-|\mu\,(z)|}\,.
\end{equation}
Given a Lebesgue measurable function $\mu:D\rightarrow {\Bbb D},$
${\Bbb D}=\{z\in {\Bbb C}: |z|<1\},$ we define the {\it maximal
dilatation} of $f$ at $z$ the function $K_{\mu}(z)$ in~(\ref{eq1}).
Note that the Jacobian of $f$ at $z\in D$ is calculated by the
formula
$$J(z,
f)=|f_z|^2-|f_{\overline{z}}|^2\,.$$
It is easy to see that
$K_{\mu_f}(z)=\frac{|f_z|+|f_{\overline{z}}|}{|f_z|-|f_{\overline{z}}|}$
whenever partial derivatives of $f$ exist at $z\in D$ and, in
addition, $J(z, f)\ne 0.$

\medskip
We will call the {\it Beltrami equation} the differential equation
of the form
\begin{equation}\label{eq2}
f_{\overline{z}}=\mu(z)\cdot f_z\,,
\end{equation}
where $\mu=\mu(z)$ is a given function with $|\mu(z)|<1$ a.a. The
{\it regular solution} of~(\ref{eq2}) in the domain $D\subset{\Bbb
C}$ is a homeomorphism $f:D\rightarrow{\Bbb C}$ of the class $W_{\rm
loc}^{1, 1}(D)$ such that $J(z, f)\ne 0$ for almost all $z\in D.$

\medskip
In the extended Euclidean space $\overline{{{\Bbb R}}^n}={{\Bbb
R}}^n\cup\{\infty\},$ we use the so-called {\it chordal metric} $h$
defined by the equalities
\begin{equation}\label{eq3C}
h(x,y)=\frac{|x-y|}{\sqrt{1+{|x|}^2} \sqrt{1+{|y|}^2}}\,,\quad x\ne
\infty\ne y\,, \quad\,h(x,\infty)=\frac{1}{\sqrt{1+{|x|}^2}}\,,
\end{equation}
see e.g.~\cite[Definition~12.1]{Va}. For a given set
$E\subset\overline{{\Bbb R}^n},$ we set
\begin{equation}\label{eq47***}
h(E):=\sup\limits_{x,y\in E}h(x, y)\,.
\end{equation}
The quantity $h(E)$ in~(\ref{eq47***}) is called the {\it chordal
diameter} of the set $E.$ As usual, the family $\frak{F}$ of
mappings $f:D\rightarrow \overline{{\Bbb C}}$ is called {\it
normal,} if from each sequence $f_n\in \frak{F},$ $n=1,2,\ldots , $
one can choose a subsequence $f_{n_k},$ $k=1,2,\ldots ,$ converging
to some mapping $f:D\rightarrow \overline{{\Bbb C}}$ locally
uniformly with respect to the metric $h.$ If, in addition, $f\in
\frak{F},$ the family $\frak{F}$ is called {\it compact}.

\medskip
A set $A\subset {\Bbb D} $ is called {\it invariantly convex} if the
set $g(A)$ is convex for any fractional-linear automorphism $g$ of
the unit disk.

Let $D$ be a domain in ${\Bbb R}^n,$ $n\geqslant 2.$ Recall some
definitions (see, for example,~\cite{KR$_1$}, \cite{KR$_2$},
\cite{IS} or \cite{ISS}). Let $\omega$ be an open set in ${\Bbb
R}^k$, $k=1,\ldots,n-1$. A continuous mapping
$\sigma\colon\omega\rightarrow{\Bbb R}^n$ is called a {\it
$k$-dimensional surface} in ${\Bbb R}^n$. A {\it surface} is an
arbitrary $(n-1)$-dimensional surface $\sigma$ in ${\Bbb R}^n.$ A
surface $\sigma$ is called {\it a Jordan surface}, if
$\sigma(x)\ne\sigma(y)$ for $x\ne y$. In the following, we will use
$\sigma$ instead of $\sigma(\omega)\subset {\Bbb R}^n,$
$\overline{\sigma}$ instead of $\overline{\sigma(\omega)}$ and
$\partial\sigma$ instead of
$\overline{\sigma(\omega)}\setminus\sigma(\omega).$ A Jordan surface
$\sigma\colon\omega\rightarrow D$ is called a {\it cut} of $D$, if
$\sigma$ separates $D,$ that is $D\setminus \sigma$ has more than
one component, $\partial\sigma\cap D=\varnothing$ and
$\partial\sigma\cap\partial D\ne\varnothing$.

A sequence of cuts $\sigma_1,\sigma_2,\ldots,\sigma_m,\ldots$ in $D$
is called {\it a chain}, if:

(i) the set $\sigma_{m+1}$ is contained in exactly one component
$d_m$ of the set $D\setminus \sigma_m,$ wherein $\sigma_{m-1}\subset
D\setminus (\sigma_m\cup d_m)$; (ii)
$\bigcap\limits_{m=1}^{\infty}\,d_m=\varnothing.$

Two chains of cuts  $\{\sigma_m\}$ and $\{\sigma_k^{\,\prime}\}$ are
called {\it equivalent}, if for each $m=1,2,\ldots$ the domain $d_m$
contains all the domains $d_k^{\,\prime},$ except for a finite
number, and for each $k=1,2,\ldots$ the domain $d_k^{\,\prime}$ also
contains all domains $d_m,$ except for a finite number.

The {\it end} of the domain $D$ is the class of equivalent chains of
cuts in $D$. Let $K$ be the end of $D$ in ${\Bbb R}^n$, then the set
$I(K)=\bigcap\limits_{m=1}\limits^{\infty}\overline{d_m}$ is called
{\it the impression of the end} $K$. Throughout the paper,
$\Gamma(E, F, D)$ denotes the family of all paths $\gamma\colon[a,
b]\rightarrow \overline{{\Bbb R}^n}$ such that $\gamma(a)\in E,$
$\gamma(b)\in F$ and $\gamma(t)\in D$ for every $t\in[a, b].$ In
what follows, $M$ denotes the modulus of a family of paths, and the
element $dm(x)$ corresponds to the Lebesgue measure in ${\Bbb R}^n,$
$n\geqslant 2,$ see~\cite{Va}. Following~\cite{Na$_2$}, we say that
the end $K$ is {\it a prime end}, if $K$ contains a chain of cuts
$\{\sigma_m\}$ such that
$\lim\limits_{m\rightarrow\infty}M(\Gamma(C, \sigma_m, D))=0$ for
some continuum $C$ in $D.$ In the following, the following notation
is used: the set of prime ends corresponding to the domain $D,$ is
denoted by $E_D,$ and the completion of the domain $D$ by its prime
ends is denoted $\overline{D}_P.$

Consider the following definition, which goes back to
N\"akki~\cite{Na$_2$}, see also~\cite{KR$_1$}--\cite{KR$_2$}. We say
that the boundary of the domain $D$ in ${\Bbb R}^n$ is {\it locally
quasiconformal}, if each point $x_0\in\partial D$ has a neighborhood
$U$ in ${\Bbb R}^n$, which can be mapped by a quasiconformal mapping
$\varphi$ onto the unit ball ${\Bbb B}^n\subset{\Bbb R}^n$ so that
$\varphi(\partial D\cap U)$ is the intersection of ${\Bbb B}^n$ with
the coordinate hyperplane.

\medskip
For the sets $A, B\subset{\Bbb R}^n$ we set, as usual,
$${\rm diam}\,A=\sup\limits_{x, y\in A}|x-y|\,,\quad {\rm dist}\,(A, B)=\inf\limits_{x\in A,
y\in B}|x-y|\,.$$
Sometimes we also write $d(A)$ instead of ${\rm diam}\,A$ and $d(A,
B)$ instead of ${\rm dist\,}(A, B),$ if no misunderstanding is
possible. The sequence of cuts $\sigma_m,$ $m=1,2,\ldots ,$ is
called {\it regular,} if
$\overline{\sigma_m}\cap\overline{\sigma_{m+1}}=\varnothing$ for
$m\in {\Bbb N}$ and, in addition, $d(\sigma_{m})\rightarrow 0$ as
$m\rightarrow\infty.$ If the end $K$ contains at least one regular
chain, then $K$ will be called {\it regular}. We say that a bounded
domain $D$ in ${\Bbb R}^n$ is {\it regular}, if $D$ can be
quasiconformally mapped to a domain with a locally quasiconformal
boundary whose closure is a compact in ${\Bbb R}^n,$ and, besides
that, every prime end in $D$ is regular. Note that space
$\overline{D}_P=D\cup E_D$ is metric, which can be demonstrated as
follows. If $g:D_0\rightarrow D$ is a quasiconformal mapping of a
domain $D_0$ with a locally quasiconformal boundary onto some domain
$D,$ then for $x, y\in \overline{D}_P$ we put:
\begin{equation}\label{eq5}
\rho(x, y):=|g^{\,-1}(x)-g^{\,-1}(y)|\,,
\end{equation}
where the element $g^{\,-1}(x),$ $x\in E_D,$ is to be understood as
some (single) boundary point of the domain $D_0.$ The specified
boundary point is unique and well-defined by~\cite[Theorem~2.1,
Remark~2.1]{IS}, cf.~\cite[Theorem~4.1]{Na$_2$}. It is easy to
verify that~$\rho$ in~(\ref{eq5}) is a metric on $\overline{D}_P,$
and that the topology on $\overline{D}_P,$ defined by such a method,
does not depend on the choice of the map $g$ with the indicated
property.

We say that a sequence $x_m\in D,$ $m=1,2,\ldots,$ converges to a
prime end of $P\in E_D$ as $m\rightarrow\infty, $ write
$x_m\rightarrow P$ as $m\rightarrow\infty,$ if for any $k\in {\Bbb
N}$ all elements $x_m$ belong to $d_k$ except for a finite number.
Here $d_k$ denotes a sequence of nested domains corresponding to the
definition of the prime end $P.$ Note that for a homeomorphism of a
domain $D$ onto $D^{\,\prime},$ the end of the domain $D$ uniquely
corresponds to some sequence of nested domains in the image under
the mapping.

\medskip
Consider the following Dirichlet problem:
\begin{equation}\label{eq2C}
f_{\overline{z}}=\mu(z)\cdot f_z\,,
\end{equation}
\begin{equation}\label{eq1A}
\lim\limits_{\zeta\rightarrow P}{\rm
Re\,}f(\zeta)=\varphi(P)\qquad\forall\,\, P\in E_D\,,
\end{equation}
where $\varphi: E_D\rightarrow {\Bbb R}$  is a predefined continuous
function. In what follows, we assume that $D$ is some simply
connected domain in ${\Bbb C}.$ The solution of the
problem~(\ref{eq2C})--(\ref{eq1A}) is called {\it regular,} if one
of two is fulfilled: or $f(z)=const$ in $D,$ or $f$ is an open
discrete $W_{\rm loc}^{1, 1}(D)$-mapping such that $J(z, f)\ne 0$
for almost any $z\in D.$

\medskip
Let $D$ be a domain in ${\Bbb R}^n,$ $n\geqslant 2.$ We say that a
function ${\varphi}:D\rightarrow{\Bbb R}$ has a {\it finite mean
oscillation} at a point $x_0\in D,$ write $\varphi\in FMO(x_0),$ if
\begin{equation}\label{eq29*!}
\limsup\limits_{\varepsilon\rightarrow
0}\frac{1}{\Omega_n\varepsilon^n}\int\limits_{B( x_0,\,\varepsilon)}
|{\varphi}(x)-\overline{{\varphi}}_{\varepsilon}|\ dm(x)<\infty\,,
\end{equation}
where
$$\overline{{\varphi}}_{\varepsilon}=\frac{1}
{\Omega_n\varepsilon^n}\int\limits_{B(x_0,\,\varepsilon)}
{\varphi}(x) dm(x)\,.$$
\medskip
Observe that, as known, $\Omega_n\varepsilon^n=m(B(x_0,
\varepsilon)),$ and that under condition (\ref{eq29*!}), a situation
is possible when
$\overline{{\varphi}_{\varepsilon}}\rightarrow\infty$ при
$\varepsilon\rightarrow 0.$
We also say that a function ${\varphi}:D\rightarrow{\Bbb R}$ has a
finite mean oscillation in $D,$ write ${\varphi}\in FMO(D),$ or
simply ${\varphi}\in FMO,$ if ${\varphi}$ has a finite mean
oscillation at any point $x_0\in D.$

\medskip
Let $M(z)\subset {\Bbb D},$ $z\in{\Bbb C}$ be some system of sets
(that is, for each $z_0\in {\Bbb C}$ the symbol $M(z_0)$ denotes
some set in ${\Bbb D}$). Denote by $\frak{M}_M$ the set of all
complex measurable functions $\mu:{\Bbb C}\rightarrow {\Bbb D},$
such that $\mu(z)\in M(z)$ for almost all $z\in {\Bbb C}.$

\medskip
We fix a point $z_0\in D$ and a function $\varphi.$ Let $M(z)\subset
{\Bbb D},$ $z\in D,$ be some system of sets. Let $\frak{F}_{\varphi,
M, z_0}(D)$ be the class of all regular solutions
$f:D\rightarrow{\Bbb C}$ of the Dirichlet problem
(\ref{eq2C})--(\ref{eq1A}), which satisfy the condition ${\rm
Im}\,f(z_0)=0$ such that $\mu\in \frak{M}_M.$ We define a function
$Q_M(z)$ by the relation
\begin{equation}\label{eq1K} Q_M(z)=\frac{1+q_M(z)}{1-q_M(z)}\,,\quad
q_M(z)=\sup\limits_{\nu\in M(z)}|\nu|\,,
\end{equation}
and we consider that $Q_M(z)\equiv 1$ for $z\in {\Bbb C} \setminus
D.$ The following statement generalizes \cite[Theorem~2]{Dyb$_1$}
for the case of arbitrary simply connected Jordan domains.

\medskip
\begin{theorem}\label{th2A}
{\sl Let $D$ be a simply connected domain in ${\Bbb C},$ and let
$\varphi$ be a continuous function in~(\ref{eq1A}). Let $M(z),$
$z\in D,$ be a family of convex compact sets, and let $Q_M$ be
integrable in $D$ and satisfies at least one of the following
conditions: either $Q_M\in FMO(\overline{D}),$ or
\begin{equation}\label{eq2G}
\int\limits_0^{\delta_0}\frac{dt}{tq_{M_{x_0}}(t)}=\infty
\end{equation}
for any $x_0\in \overline{D}$ and some $\delta_0=\delta(x_0)>0,$
where
$q_{M_{x_0}}(t)=\frac{1}{2\pi}\int\limits_0^{2\pi}Q_M(x_0+te^{i\theta})\,d\theta.$
Then the family $\frak{F}_{\varphi, M, z_0}(D)$ is compact in $D.$ }
\end{theorem}

\section{Preliminaries}

In what follows, $M(\Gamma)$ denotes the conformal modulus of a
family of paths $\Gamma$ (see \cite[section~6]{Va}), and the element
$dm(x)$ corresponds to a Lebesgue measure in ${\Bbb R}^n,$
$n\geqslant 2,$ see~\cite{Va}. Everywhere below, unless otherwise
stated, the boundary and the closure of a set are understood in the
sense of an extended Euclidean space $\overline{{\Bbb R}^n}.$ Let
$x_0\in\overline{D},$ $x_0\ne\infty,$
$$S(x_0,r) = \{
x\,\in\,{\Bbb R}^n : |x-x_0|=r\}\,, S_i=S(x_0, r_i)\,,\quad
i=1,2\,,$$
\begin{equation}\label{eq6} A=A(x_0, r_1, r_2)=\{ x\,\in\,{\Bbb R}^n :
r_1<|x-x_0|<r_2\}\,.
\end{equation}
Let $Q:{\Bbb R}^n\rightarrow {\Bbb R}^n$ be a Lebesgue measurable
function satisfying the condition $Q(x)\equiv 0$ for $x\in{\Bbb
R}^n\setminus D,$ and let $p\geqslant 1.$ Given sets $E$ and $F$ and
a given domain $D$ in $\overline{{\Bbb R}^n}={\Bbb R}^n\cup
\{\infty\},$ we denote by $\Gamma(E, F, D)$ the family of all paths
$\gamma:[0, 1]\rightarrow \overline{{\Bbb R}^n}$ joining $E$ and $F$
in $D,$ that is, $\gamma(0)\in E,$ $\gamma(1)\in F$ and
$\gamma(t)\in D$ for all $t\in (0, 1).$ According
to~\cite[Сhap.~7.6]{MRSY}, a mapping $f:D\rightarrow \overline{{\Bbb
R}^n}$ is called a {\it ring $Q$-mapping at the point $x_0\in
\overline{D}\setminus \{\infty\}$}, if the condition
\begin{equation} \label{eq3*!gl0}
M(f(\Gamma(S_1, S_2, D)))\leqslant \int\limits_{A\cap D} Q(x)\cdot
\eta^n (|x-x_0|)\, dm(x)
\end{equation}
holds for all $0<r_1<r_2<d_0:=\sup\limits_{x\in D}|x-x_0|$ and all
Lebesgue measurable functions $\eta:(r_1, r_2)\rightarrow [0,
\infty]$ such that
\begin{equation}\label{eq*3gl0}
\int\limits_{r_1}^{r_2}\eta(r)\,dr\geqslant 1\,.
\end{equation}
The mapping $f:D\rightarrow \overline{{\Bbb R}^n}$ is called a {\it
ring $Q$-mapping in $\overline{D}\setminus\{\infty\}$ }
if~(\ref{eq3*!gl0}) holds for any~$x_0\in
\overline{D}\setminus\{\infty\}.$ This definition can also be
applied to the point~$x_0=\infty$ by inversion:
$\varphi(x)=\frac{x}{|x|^2},$ $\infty\mapsto 0.$ In what follows,
$h$ denotes the so-called chordal metric defined by~(\ref{eq3C}).

\medskip
The next important lemma follows by~\cite[Theorems~4.1 and
4.2]{RSS}.

\medskip\begin{lemma}\label{lem1}
{\sl \, Let $D$ be a domain in ${\Bbb R}^n,$ $n\geqslant 2,$ and let
$Q:D\rightarrow [1, \infty]$ be a Lebesgue measurable function. In
addition, let, $f_k,$ $k=1,2, \ldots$ be a sequence of
homeomorphisms of $D$ into ${\Bbb R}^n,$ which satisfy
conditions~(\ref{eq3*!gl0})--(\ref{eq*3gl0} at any point $x_0\in D$
that converges to some mapping $f:D\rightarrow \overline{{\Bbb
R}^n}$ locally uniformly in $D$ with respect to the chordal metric
$h.$ Assume that the function $Q$ satisfies at least one of two
following conditions: either  $Q\in FMO(D),$ or
\begin{equation}\label{eq2H}
\int\limits_0^{\delta_0}\frac{dt}{tq^{\frac{1}{n-1}}_{x_0}(t)}=\infty\,,
\end{equation}
for any $x_0\in D$ and some $\delta_0=\delta(x_0)>0,$ where
$q_{x_0}(t)=\frac{1}{\omega_{n-1}t^{n-1}}\int\limits_{S(x_0,
t)}Q(x)\,d\mathcal{H}^{n-1},$ and $\mathcal{H}^{n-1}$ denotes
$(n-1)$-dimensional Hausdorff measure. Then the mapping $f$ is
either a homeomorphism $f:D\rightarrow {\Bbb R}^n,$ or a constant $c
\in \overline{{\Bbb R}^n}.$}\end{lemma}

\medskip
Let $I$ be a fixed set of indices and let $D_i,$ $i\in I,$ be some
sequence of domains. Following~\cite[Sect.~2.4]{NP}, we say that a
family of domains $\{D_i\}_{i\in I}$ is {\it equi-uniform} if for
any $r> 0$ there exists a number $\delta> 0$ such that the
inequality
\begin{equation}\label{eq17***}
M(\Gamma(F^{\,*},F, D_i))\geqslant \delta
\end{equation}
holds for any $i\in I$ and any continua $F, F^*\subset D$ such that
$h(F)\geqslant r$ and $h(F^{\,*})\geqslant r.$ If $D$ is one domain
satisfying condition~(\ref{eq17***}), then it is called {\it
uniform}.

\medskip
Given numbers $\delta>0,$ a domain $D\subset {\Bbb R}^n,$
$n\geqslant 2,$ a point $a\in D$ and a Lebesgue measurable function
$Q:{\Bbb R}^n\rightarrow {\Bbb R}^n,$ $Q(x)\equiv 0$ for $x\in{\Bbb
R}^n\setminus D,$ denote by $\frak{F}_{Q, a, \delta}(D)$ the family
of all homeomorphisms $f:D\rightarrow \overline{{\Bbb R}^n}$
satisfying (\ref{eq3*!gl0})--(\ref{eq*3gl0}) in $\overline{D}$ such
that $h(f(a),
\partial f(D))\geqslant\delta,$ $h(\overline{{\Bbb R}^n}\setminus
f(D))\geqslant \delta.$ The following statement holds
(see~\cite[Theorem~2]{SevSkv}).

\medskip
\begin{lemma}\label{th1} {\sl\, Let $D$ be regular,
and let $D_f^{\,\prime}=f(D)$ be bounded domains with a locally
quasiconformal boundary which are equi-uniform over all
$f\in\frak{F}_{Q, a, \delta}(D).$ If $Q\in FMO(D),$ or the
condition~(\ref{eq2H}) holds, then any $f\in\frak{F}_{Q, a,
\delta}(D)$ has a continuous extension~$\overline{f}:
\overline{D}_P\rightarrow \overline{{\Bbb R}^n}$ and, in addition,
the family $f\in\frak{F}_{Q, a, \delta}(\overline{D})$ of all
extended mappings $\overline{f}: \overline{D}_P\rightarrow
\overline{{\Bbb R}^n}$ is equicontinuous in $\overline{D}_P.$
  }
\end{lemma}

\medskip
We also need to formulate a similar statement for homeomorphisms
inverse to~(\ref{eq3*!gl0}). For this purpose, consider the
following definition.

\medskip
For domains $D\subset {\Bbb R}^n$ and
$D^{\,\prime}\subset\overline{{\Bbb R}^n},$ $n\geqslant 2,$ points
$a\in D,$ $b\in D^{\,\prime}$ and a Lebesgue measurable function $Q:
{\Bbb R}^n\rightarrow [0, \infty]$ such that $Q(x)\equiv 0$ for
$x\not\in D,$ we denote by ${\frak S}_{a, b, Q}(D, D^{\,\prime})$
the family of all homeomorphisms $h$ of $D^{\,\prime}$ onto $D$ such
that the mapping $f=h^{\,-1}$ satisfies the
condition~(\ref{eq3*!gl0}) in $\overline{D},$ while $f(a)=b.$

\medskip
The boundary of the domain $D$ is called {\it weakly flat at the
point $x_0,$} if for every number $P>0$ and every neighborhood $U$
of this point there is a neighborhood $V$ of point $x_0$ such that
$M(\Gamma(E, F, D))>P$ for any continua $E$ and $F,$ satisfying
conditions $F\cap \partial U\ne\varnothing\ne F\cap
\partial V.$ The boundary of domain
$D$ is called {\it weakly flat} if it is such at each of its point.
The following statement holds (see e.g.~\cite[Theorem~7.1]{SSD}).

\medskip
\begin{lemma}\label{th4}
{\sl Assume that $D$ is a regular domain, and that $D^{\,\prime}$
has a weakly flat boundary, none of the components of which
degenerates into a point. Then any $h\in {\frak S}_{a, b, Q}(D,
D^{\,\prime})$ has a continuous extension
$\overline{h}:\overline{D^{\,\prime}}\rightarrow \overline{D}_P,$
while $\overline{h}(\overline{D^{\,\prime}})=\overline{D}_P$ and, in
addition, the family ${\frak S}_{a, b, Q}(\overline{D},
\overline{D^{\,\prime}})$ of all extended mappings
$\overline{h}:\overline{D^{\,\prime}}\rightarrow \overline{D}_P$ is
equicontinuous in $\overline{D^{\,\prime}}.$ }
\end{lemma}

\section{Proof of Theorem~\ref{th2A}}

In general, we will use the scheme of proving Theorem~1.2
in~\cite{SD}.

\medskip
{\textbf I.} Let $f_m\in \frak{F}^{\mathcal{M}}_{\varphi, \Phi,
z_0}(D),$ $m=1,2,\ldots .$ By Stoilow's factorization theorem (see,
e.g., \cite[5(III).V]{St}) a mapping $f_m$ has a representation
\begin{equation}\label{eq2E}
f_m=\varphi_m\circ g_m\,,
\end{equation}
where $g_m$ is some homeomorphism, and $\varphi_m$ is some analytic
function. By Lemma~1 in~\cite{Sev}, the mapping $g_m$ belongs to the
Sobolev class $W_{\rm loc}^{1, 1}(D)$ and has a finite distortion.
Moreover, by~\cite[(1).C, Ch.~I]{A}
\begin{equation}\label{eq1B}
{f_m}_z={\varphi_m}_z(g_m(z)){g_m}_z,\qquad
{f_m}_{{\overline{z}}}={\varphi_m}_z(g_m(z)){g_m}_{\overline{z}}
\end{equation}
for almost all $z\in D.$ Therefore, by the relation~(\ref{eq1B}),
$J(z, g_m)\ne 0$ for almost all $z\in D,$ in addition,
$K_{\mu_{f_m}}(z)=K_{\mu_{g_m}}(z).$

\medskip
\textbf{II.} We prove that $\partial g_m (D)$ contains at least two
points. Suppose the contrary. Then either $g_m(D)={\Bbb C},$ or
$g_m(D)={\Bbb C}\setminus \{a\},$ where $a\in {\Bbb C}.$ Consider
first the case $g_m(D)={\Bbb C}.$ By Picard's theorem
$\varphi_m(g_m(D))$ is the whole plane, except perhaps one point
$\omega_0\in {\Bbb C}.$ On the other hand, for every $m=1,2,\ldots$
the function $u_m(z):={\rm Re}\,f_m(z)={\rm
Re}\,(\varphi_m(g_m(z)))$ is continuous on the compact set
$\overline{D}$ under the condition~(\ref{eq1A}) by the continuity
of~$\varphi.$ Therefore, there exists $C_m>0$ such that $|{\rm
Re}\,f_m(z)|\leqslant C_m$ for any $z\in D,$ but this contradicts
the fact that $\varphi_m(g_m(D))$ contains all points of the complex
plane except, perhaps, one. The situation $g_m(D)={\Bbb C}\setminus
\{a\},$ $a\in {\Bbb C},$ is also impossible, since the domain
$g_m(D)$ must be simply connected in ${\Bbb C}$ as a homeomorphic
image of the simply connected domain $D.$

\medskip
Therefore, the boundary of the domain $g_m(D)$ contains at least two
points. Then, according to Riemann's mapping theorem, we may
transform the domain $g_m(D)$ onto the unit disk ${\Bbb D}$ using
the conformal mapping $\psi_m.$ Let $z_0\in D $ be a point from the
condition of the theorem. By using an auxiliary conformal mapping
$$\widetilde{\psi_m}(z)=\frac{z-(\psi_m\circ
g_m)(z_0)}{1-z\overline{(\psi_m\circ g_m)(z_0)}}$$ of the unit disk
onto itself we may consider that $(\psi_m\circ g_m)(z_0)=0.$ Now,
by~(\ref{eq2E}) we obtain that
$$
f_m=\varphi_m\circ g_m= \varphi_m\circ\psi^{\,-1}_m\circ\psi_m\circ
g_m=F_m\circ G_m\,,\quad m=1,2,\ldots\,,
$$
where $F_m:=\varphi_m\circ\psi^{\,-1}_m,$ $F_m:{\Bbb D}\rightarrow
{\Bbb C},$ and $G_m=\psi_m\circ g_m.$
Obviously, a function $F_m$ is analytic, and $G_m$ is a regular
Sobolev homeomorphism in $D.$ In particular, ${\rm Im}\,F_m(0)=0$
for any $m\in {\Bbb N}.$

\medskip
\textbf{III.} Observe that
\begin{equation}\label{eq7A}
\int\limits_D K_{\mu_{G_m}}(z)\,dm(z) \leqslant \int\limits_D
Q_M(z)\,dm(z) <\infty\,,
\end{equation}
because by the condition $\mu(z)\in M(z)$ for almost any $z\in D,$
moreover, the inequality $K_{\mu_{G_m}}(z)\leqslant Q_M(z)$ holds
for almost any $z\in D,$ where $Q_M(z)$ does not depend on
$m=1,2,\ldots $ and is integrable.

\medskip
\textbf{IV.} We prove that each map $G_m,$ $m=1,2,\ldots ,$ has a
continuous extension to $E_D,$ in addition, the family of extended
maps $\overline{G}_m,$ $m=1,2,\ldots ,$ is equicontinuous in
$\overline{D}_P.$ Indeed, as proved in item~\textbf{III},
$K_{\mu_{G_m}}\in L^1(D).$ By~\cite[Theorem~3]{KPRS} (see
also~\cite[Theorem~3.1]{LSS}) each $G_m,$ $m=1,2,\ldots, $ is a ring
$Q$-homeomorphism in $\overline{D}$ for $Q=K_{\mu_{G_m}}(z),$ where
$\mu$ is defined in~(\ref{eq2C}), and $K_{\mu}$ my be calculated by
the formula~(\ref{eq1}). Note that the unit disk ${\Bbb D}$ is a
uniform domain as a finitely connected flat domain at its boundary
with a finite number of boundary components (see, for example,
~\cite[Theorem~6.2 and Corollary~6.8]{Na$_1$}). Then the desirable
the conclusion is a statement of Lemma~\ref{th1}.

\medskip
\textbf{V.} Observe that, the inverse homeomorphisms $G^{\,-1}_m,$
$m=1,2,\ldots ,$ have a continuous extension $\overline{G}^{\,-1}_m$
to $\partial {\Bbb D}$ in terms of prime ends in $D,$ and
$\{\overline{G}_m^{\,-1}\}_{m=1}^{\infty}$ is equicontinuous in
$\overline{\Bbb D}$ as a family of mappings from $\overline{\Bbb D}$
to $\overline{D}_P.$ Indeed, by the item~\textbf{IV} mappings $G_m,$
$m=1,2,\ldots, $ are ring $K_{\mu_{G_m}}(z)$-homeomorphisms in $D$
such that $G^{\,-1}_m(0)=z_0$ for any $m=1,2,\ldots .$ In this case,
the possibility of a continuous extension of $G^{\,-1}_m$ to
$\partial {\Bbb D},$ and the equicontinuity of
${\{\overline{G}_m^{\,-1}\}}_{m=1}^{\infty}$ as mappings
$G_m^{\,-1}:\overline{\Bbb D}\rightarrow \overline{D}_P$ follows by
Lemma~\ref{th4}.

\medskip
\textbf{VI.} Since, as proved above the family
${\{G_m\}}^{\infty}_{m=1}$ is equicontinuous in~$D,$ according to
Arzela-Ascoli criterion there exists an increasing subsequence of
numbers $m_k,$ $k=1,2,\ldots ,$ such that $G_{m_k}$ converges
locally uniformly in $D$ to some continuous mapping $G:D\rightarrow
\overline{{\Bbb C}}$ as $k\rightarrow\infty$  (see, e.g.,
\cite[Theorem~20.4]{Va}). By Lemma~\ref{lem1}, either $G$ is a
homeomorphism with values in ${\Bbb R}^n,$ or a constant
in~$\overline{{\Bbb R}^n}.$ Let us prove that the second case is
impossible. Let us apply the approach used in proof of the second
part of Theorem~21.9 in~\cite{Va}. Suppose the contrary: let
$G_{m_k}(x)\rightarrow c=const$ as $k\rightarrow\infty.$ Since
$G_{m_k}(z_0)=0$ for all $k=1,2,\ldots ,$ we have that $c=0.$ By
item~\textbf{V}, the family of mappings~$G^{\,-1}_m,$ $m=1,2,\ldots
,$ is equicontinuous in ${\Bbb D}.$ Then
$$h(z, G^{\,-1}_{m_k}(0))=h(G^{\,-1}_{m_k}(G_{m_k}(z)), G^{\,-1}_{m_k}(0))\rightarrow 0$$
as $k\rightarrow\infty,$ which is impossible because $z$ is an
arbitrary point of the domain $D.$ The obtained contradiction
refutes the assumption made above. Thus, $G:D\rightarrow {\Bbb C}$
is a homeomorphism.

\medskip
\textbf{VII.} According to~\textbf{V}, the family of mappings
$\{\overline{G}_m^{\,-1}\}_{m=1}^{\infty}$ is equicontinuous
in~$\overline{\Bbb D}.$ By the Arzela-Ascoli criterion (see, e.g.,
\cite[Theorem~20.4]{Va}) we may consider that
$\overline{G}^{\,-1}_{m_k}(y),$ $k=1,2,\ldots, $ converges to some
mapping $\widetilde{F}:\overline{{\Bbb D}}\rightarrow \overline{D}$
as $k\rightarrow\infty$ uniformly in~$\overline{D}.$ Let us to prove
that $\widetilde{F}=\overline{G}^{\,-1}.$ For this purpose, we show
that~$G(D)={\Bbb D}.$ Fix $y\in {\Bbb D}.$ Since $G_{m_k}(D)={\Bbb
D}$ for every $k=1,2,\ldots, $ we obtain that $G_{m_k}(x_k)=y$ for
some $x_k\in D.$ Since $D$ is regular, the metric space
$(\overline{D}_P, \rho)$ is compact. Thus, we may assume that
$\rho(x_k, x_0)\rightarrow 0$ as $k\rightarrow\infty,$ where
$x_0\in\overline{D}_P.$ By the triangle inequality and the
equicontinuity of ${\{\overline{G}_m\}}^{\infty}_{m=1}$ onto
$\overline{D}_P,$ see~\textbf{IV}, we obtain that
$$|\overline{G}
(x_0)-y|=|\overline{G}(x_0)-\overline{G}_{m_k}(x_k)|\leqslant
|\overline{G}(x_0)-\overline{G}_{m_k}(x_0)|+|\overline{G}_{m_k}(x_0)-\overline{G}_{m_k}(x_k)|\rightarrow
0$$
as $k\rightarrow\infty.$ Thus, $\overline{G}(x_0)=y.$ Observe that
$x_0\in D,$ because $G$ is a homeomorphism. Since $y\in {\Bbb D}$ is
arbitrary, the equality $G(D)={\Bbb D}$ is proved. In this case,
$G^{\,-1}_{m_k}\rightarrow G^{\,-1}$ locally uniformly in ${\Bbb D}$
as $k\rightarrow\infty$ (see, e.g., \cite[Lemma~3.1]{RSS}). Thus,
$\widetilde{F}(y)=G^{\,-1}(y)$ for every $y\in {\Bbb D}.$

Finally, since $\widetilde{F}(y)=G^{\,-1}(y)$ for any $y\in {\Bbb
D}$ and, in addition, $\widetilde{F}$ has a continuous extension on
$\partial{\Bbb D},$ due to the uniqueness of the limit at the
boundary points we obtain that
$\widetilde{F}(y)=\overline{G}^{\,-1}(y)$ for $y\in \overline{{\Bbb
D}}.$ Therefore, we have proved
that~$\overline{G}^{\,-1}_{m_k}\rightarrow \overline{G}^{\,-1}$
uniformly in~$\overline{\Bbb D}$ with as $k\rightarrow\infty$ with
respect to the metrics $\rho$ in $\overline{D}_P.$

\medskip
\textbf{VIII.} By~\textbf{VII,} for $y=e^{i\theta}\in
\partial {\Bbb D}$
\begin{equation}\label{eq4E}
{\rm
Re\,}F_{m_k}(e^{i\theta})=\varphi\left(\overline{G}^{\,-1}_{m_k}(e^{i\theta})\right)\rightarrow
\varphi\left(\overline{G}^{\,-1}(e^{i\theta})\right)
\end{equation}
as $k\rightarrow\infty$ uniformly on $\theta\in [0, 2\pi).$ Since by
the construction ${\rm Im\,}F_{m_k}(0)=0$ for any $k=1,2,\ldots,$ by
the Schwartz formula (see, e.g., \cite[section~8.III.3]{GK}) the
analytic function $F_{m_k}$ is uniquely restored by its real part,
namely,
\begin{equation}\label{eq4A}
F_{m_k}(y)=\frac{1}{2\pi i}\int\limits_{S(0,
1)}\varphi\left(\overline{G}^{\,-1}_{m_k}(t)\right)\frac{t+y}{t-y}\cdot\frac{dt}{t}\,.
\end{equation}
Set
\begin{equation}\label{eq5B}
F(y):=\frac{1}{2\pi i}\int\limits_{S(0,
1)}\varphi\left(\overline{G}^{\,-1}(t)\right)\frac{t+y}{t-y}\cdot\frac{dt}{t}\,.
\end{equation}
Let $K\subset {\Bbb D}$ be an arbitrary compact set, and let $y\in
K.$ By~(\ref{eq4A}) and~(\ref{eq5B}) we obtain that
\begin{equation}\label{eq11A}
|F_{m_k}(y)-F(y)|\leqslant \frac{1}{2\pi}\int\limits_{S(0,
1)}\bigl|\varphi(\overline{G}^{\,-1}_{m_k}(t))-\varphi(\overline{G}
^{\,-1}(t))\bigr|\left|\frac{t+y}{t-y}\right|\,|dt|\,.
\end{equation}
Since $K$ is compact, there is $0<R_0=R_0(K)<\infty$ such that
$K\subset B(0, R_0).$ By the triangle inequality $|t+y|\leqslant
1+R_0$ and $|t-y|\geqslant |t|-|y|\geqslant 1-R_0$ for $y\in K$ and
any $t\in {\Bbb S}^1.$ Thus
$$
\left|\frac{t+y}{t-y}\right|\leqslant \frac{1+R_0}{1-R_0}:=M=M(K)\,.
$$
Put $\varepsilon>0.$ By~(\ref{eq4E}), for a number
$\varepsilon^{\,\prime}:=\frac{\varepsilon}{M}$ there is
$N=N(\varepsilon, K)\in {\Bbb N}$ such that
$\bigl|\varphi\left(\overline{G}^{\,-1}_{m_k}(t)\right)-\varphi\left
(\overline{G}^{\,-1}(t)\right)\bigr|<\varepsilon^{\,\prime}$ for any
$k\geqslant N(\varepsilon)$ and $t\in {\Bbb S}^1.$ Now,
by~(\ref{eq11A})
\begin{equation}\label{eq13A}
|F_{m_k}(y)-F(y)|<\varepsilon \quad \forall\,\,k\geqslant N\,.
\end{equation}
It follows from~(\ref{eq13A}) that the sequence $F_{m_k}$ converges
to $F$ as $k\rightarrow\infty$ in the unit disk locally uniformly.
In particular, we obtain that ${\rm Im\,}F(0)=0.$ Note that $F$ is
analytic function in ${\Bbb D}$ (see remarks made at the end of
item~8.III in~\cite{GK}), and

$${\rm Re}\,F(re^{i\psi})=\frac{1}{2\pi}\int\limits_0^{2\pi}
\varphi\left(\overline{G}^{\,-1}(e^{i\theta})\right)\frac{1-r^2}{1-2r\cos(\theta-\psi)+r^2}\,d\theta$$
for $z=re^{i\psi}.$
By~\cite[Theorem~2.10.III.3]{GK}
\begin{equation}\label{eq15A}
\lim\limits_{\zeta\rightarrow z}{\rm
Re}\,F(\zeta)=\varphi(\overline{G}^{\,-1}(z))\quad\forall\,\,z\in\partial
{\Bbb D}\,.
\end{equation}
Observe that $F$ either is a constant or open and discrete (see,
e.g., \cite[Ch.~V,  I.6 and II.5]{St}). Thus, $f_{m_k}=F_{m_k}\circ
G_{m_k}$ converges to $f=F\circ G$ locally uniformly as
$k\rightarrow\infty,$ where $f=F\circ G$ either is a constant or
open and discrete. Moreover, by~(\ref{eq15A})
$$\lim\limits_{\zeta\rightarrow P}{\rm Re\,}
f(\zeta)= \lim\limits_{\zeta\rightarrow P}{\rm
Re\,}F(G(\zeta))=\varphi(G^{\,-1}(G(P)))=\varphi(P)\,.$$

\textbf{IX.}  Since by~\textbf{VI} $G$ is a homeomorphism,
by~\cite[Lemma~1 and Theorem~1]{L$_2$} $G$  is a regular solution of
the equation~(\ref{eq2C}) for some function~$\mu:{\Bbb C}\rightarrow
{\Bbb D}.$ Since the set of points of the function $F,$ where its
Jacobian is zero, consist only of isolated points (see~\cite[Ch.~V,
5.II and 6.II]{St}), $f$ is regular solution of the Dirichlet
problem~(\ref{eq2C})--(\ref{eq1A}) whenever $F\not\equiv const.$ It
remains to show that $\mu\in \frak{M}_M.$ If $f(z)=c=const$ in $D,$
due to the condition~(\ref{eq1A}) we obtain that $f_n(z)=c$ in $D,$
and $\mu_n(z)=0\in M(z)$ for almost any $z\in D.$ In this case,
$\mu(z)=0$ for almost any $z\in D,$ in particular, $\mu\in
\frak{M}_M.$

\medskip
Let $f(z)\ne const.$ As proved above, $f$ is regular. Since $f_n(z)$
converge to $f(z)$ locally uniformly in $D,$ in addition, the
jacobian of $f$ does not vanish almost everywhere,
by~\cite[Lemma~1]{L$_2$} $\mu(z)\in {\rm inv\,co} M_0(z)$ for almost
any $z\in D,$ where ${\rm inv\,co}\,A$ denotes the invariant-convex
hull of the set $A\subset {\Bbb C}$ (see, e.g., \cite{Ryaz}), and
$M_0(z)$ is a cluster set of $\mu_n(z),$ $n=1,2,\ldots .$ Obviously,
there is a set $D_0\subset D$ such that $\mu_n(z)\in M(z)$ and
$\mu(z)\in {\rm inv\,co}\, M_0(z)$ for all $z\in D_0$ and any $n\in
{\Bbb N},$ where $m(D\setminus D_0)=0.$ Fix $z_0\in D_0.$ Let
$w_0\in M_0(z_0).$ Then there is a subsequence of numbers $n_{k},$
$k=1,2,\ldots ,$ for which $\mu_{n_k}(z_0)$ converge as
$k\rightarrow\infty$ and $\lim\limits_{k\rightarrow
\infty}\mu_{n_k}(z_0)=w_0.$ Since, by the assumption,
$\mu_{n_k}(z_0)\in M(z_0)$ for any $k=1,2,\ldots ,$ in addition,
$M(z_0)$ is closed, we obtain that $w_0\in M(z_0).$ Thus,
\begin{equation}\label{eq1D}
M_0(z_0)\subset M(z_0)\,.
\end{equation}
It follows from~(\ref{eq1D}) that
\begin{equation}\label{eq2J}
{\rm inv\,co}\,M_0(z_0)\subset M(z_0)\,,
\end{equation}
because $M(z_0)$ is invariant-convex. Thus,
$$\mu(z_0)\in {\rm inv\,co}\,M_0(z_0)\subset M(z_0)$$ for almost any
$z_0\in D,$ that is desired conclusion.~$\Box$

\medskip
\medskip

\medskip
{\bf \noindent Oleksandr Dovhopiatyi} \\
{\bf 1.} Zhytomyr Ivan Franko State University,  \\
40 Bol'shaya Berdichevskaya Str., 10 008  Zhytomyr, UKRAINE \\
alexdov1111111@gmail.com

\medskip
{\bf \noindent Evgeny Sevost'yanov} \\
{\bf 1.} Zhytomyr Ivan Franko State University,  \\
40 Bol'shaya Berdichevskaya Str., 10 008  Zhytomyr, UKRAINE \\
{\bf 2.} Institute of Applied Mathematics and Mechanics\\
of NAS of Ukraine, \\
1 Dobrovol'skogo Str., 84 100 Slavyansk,  UKRAINE\\
esevostyanov2009@gmail.com

\end{document}